\documentclass[preprint,12pt,authoryear]{elsarticle}

\usepackage{graphicx} 
\usepackage{amsmath} 
\usepackage{amssymb}  
\usepackage{subfigure}
\newtheorem{definition}{\textbf{Definition}}
\newtheorem{remark}{\textbf{Remark}}

\newtheorem{assumption}{\textbf{Assumption}}
\newtheorem{example}{\textbf{Example}}
\newtheorem{proposition}{\textbf{Proposition}}
\usepackage[overload]{empheq} 

\usepackage{natbib}

\usepackage{algorithm,algorithmic}

\begin{document}

\begin{frontmatter}

\title{An approach to parameter identifiability for a class of nonlinear models represented in LPV form}


\author{Krishnan Srinivasarengan, Jos\'{e} Ragot, Christophe Aubrun and Didier Maquin}

\address{Universit\'{e} de Lorraine, CNRS, CRAN, F-54500, Nancy, France, {\tt\small didier.maquin@univ-lorraine.fr}}

\begin{abstract}
In several model-based system maintenance problems, parameters are used to represent unknown characteristics of a component, equipment degradation, etc. This allows for modelling constant, slow-varying terms. The identifiability of these parameters is an important condition   to estimate them. Linear Parameter Varying (LPV) models are being increasingly used in the industries as a bridge between linear and nonlinear models. Techniques exist that can rewrite some nonlinear models in LPV form.  However, the problem of identifiability of these models is still at a nascent stage. In this paper, we propose an approach to verify identifiability of unknown parameters for LPV state-space models. It makes use of a parity-space like formulation to eliminate the states of the model. The resulting input-output-parameter equation is analysed to verify the identifiability of the original model or a subset of unknown parameters. This approach provides a framework for both continuous-time and discrete-time models and we illustrate it using examples.  
\end{abstract}


\end{frontmatter}

\section{Introduction}
Parameter identifiability studies are motivated by the need for well-posed problems in several applications. For instance, to estimate a parameter through optimization (see \cite{beelen2017joint}), we need a unique set of parameters to satisfy the inputs and outputs of the experiment (at least locally). It is in this respect that the distinguishability property of parameters is defined and forms the basis of the parameter identifiability analysis. Consider the general nonlinear system model of the form,
\begin{subequations}
\label{eq_ch_id_gen_nonlinear_model}
  \begin{align}[left ={\Sigma_\theta : \empheqlbrace}]
    \dot{x}(t) &= f(x(t), u(t), \theta) \label{eq_ch_id_gen_nonlinear_model_state}\\
    y(t) &= h(x(t), u(t), \theta)\label{eq_ch_id_gen_nonlinear_model_output}
   \end{align}
\end{subequations}%
with $x \in \mathbb{R}^n$, $u \in \mathbb{R}^m$, $y \in \mathbb{R}^p$ and the constant parameters $\theta \in \mathbb{R}^q$. Distinguishability property of a model structure refers to (\cite{ljung1994global})
\begin{align}
\tilde{y}(t|\theta') \equiv \tilde{y}(t|\theta'') \Rightarrow \theta' = \theta'' \label{eq_ch_id_distinguishability}
\end{align}%
where, $\tilde{y}$ stands for the output \eqref{eq_ch_id_gen_nonlinear_model_output} computed as the solution of the system \eqref{eq_ch_id_gen_nonlinear_model} for an input $\tilde{u}$ and $\theta'$ (or ${\theta''}$) as the parameter. The essence of distinguishability is captured by the property of parameter identifiability, which refers to whether the model parameter(s) can be uniquely identified by a set of input-output data.  It is to be noted that parameter identifiability property assumes error-free model and noise-free data and hence is not a sufficient condition for the existence of a solution.

To understand the relevance of this property for applications such as fault diagnosis or prognosis, consider the problem of equipment degradation estimation in large systems. Maintenance activities are rare and typically involve gathering data by deploying temporary sensors (e.g., hand-held monitors). Consequently, only a finite amount of data is available to estimate relevant parameters that are indicative of the underlying degradation phenomena. If a model-based estimation is used, the procedure requires that the available input-output data admit a unique solution for the set of unknown parameters, at least within a known range of parameter values. This constraint can be reformulated as  local parameter identifiability of the model. 

To start with, we review the literature on the available methods to verify identifiability, focusing on nonlinear models and those amenable to LPV forms, in the Sec.~\ref{sec_relevant_lit}. This review includes the limited, but relevant works for the identifiability of discrete-time models as well. This is followed by formulation of the problem tackled in this paper in Sec.~\ref{sec_problem_formulation}. The literature contains several definitions for identifiability and this section outlines those definitions of interest, their relationships as well as the assumptions that underlie the identifiability results in this paper. The Sec.~\ref{sec_id_algo_continuous} presents the proposed algorithm in a step-by-step manner for continuous-time models with examples illustrating the procedure. The discrete-time models are treated similarly, with the focus on the steps different from that of continuous-time models, in Sec.~\ref{sec_id_algo_discrete}. A discussion on the systematic implementation of the proposed approach is given in Sec.~\ref{sec_ch_id_algo} followed by some concluding remarks in Sec.~\ref{sec_conclusions}.

\section{Relevant literature} \label{sec_relevant_lit}
To characterize identifiability of a nonlinear model encompassed by \eqref{eq_ch_id_gen_nonlinear_model}, there are several definitions in the literature. These definitions and the corresponding characterizations vary due to several factors, including,
\begin{itemize}
	\item The characteristics of the functions $f$ and $h$ (e.g., analytic, homogeneous, meromorphic)
	\item The characteristics of the inputs (sufficiently continuous/differentiable, piecewise continuous etc.)
	\item Assumptions on the initial (state) conditions 
	\item Neighbourhood of the identifiability characterization (local around a particular $\theta$ or global)
\end{itemize}
There are also nuances associated with strong and weak notions of the identifiability. For example, in \cite{nemcova_structural_2010}, the authors note that their identifiability definitions are weaker compared to that in \cite{xia2003identifiability} because the distinguishability in the latter work is for `all the admissible inputs within an open dense subset', whereas it is just `at least one input' as a sufficient condition in their own work\footnote{But at the same time, the authors in \cite{nemcova_structural_2010} consider piece-wise continuous inputs, which are more representative in biological system identification.}. 

\subsection{Identifiability of continuous-time models} 
Starting with the study of structural identifiabiality of linear models in \cite{bellman_structural_1970}, several methods have explored the problem in the following decades. In a broad sense, the methods to verify identifiability could be classified as those that perform:
\begin{enumerate}
	\item Analysis of observables
	\item Analysis of the system map
\end{enumerate}
The term \emph{observables} has been borrowed from \cite{chis_structural_2011} and roughly refers to the outputs and the parameter information embedded in them. That is, the first classification refers to verifying directly, whether the outputs (and inputs) provide a way to validate the distinguishability property in \eqref{eq_ch_id_distinguishability}. Methods such as Taylor series approach and generating series approaches fall into this category. The second class of methods look at some specific properties of the system model to check for identifiability. Isomorphism based approaches or approaches that consider identifiability as an extended observability property belong to the second class. This classification is not strict as several methods cross over. For example, the implicit function theorem based approach (\cite{xia2003identifiability}) and the differential algebraic tools based approach (\cite{bellu2007daisy}) exploit the system model properties to eliminate the latent (state) variables and then analyze the observables.

\subsubsection{Taylor series and generating series approach}
The Taylor series approach is one of the first proposed for identifiability analysis in \cite{pohjanpalo1978system}. By considering the system output as an analytic function of time, it exploits the fact that their derivatives should hold all possible information about the unknown parameters $\theta$. The uniqueness of the Taylor series expansion of this function is an indication of the identifiability of the system. If the test fails, more coefficients are to be computed and verified again.

The generating series approach in \cite{walter1982global} is conceptually similar to the Taylor series approach and is applicable to control affine models of the form,
\begin{align}
\dot{x}(t) &= f(x(t),\theta) + g(x(t),\theta) u(t), \qquad x(t_0) = x_0(\theta)  \nonumber \\
y(t) &= h(x(t),\theta) \label{eq_ch_id_control_affine}
\end{align}
Instead of using derivatives, Lie derivative expansion of the output functions along  the vector fields $f$ and $g$ are computed. The coefficients of the output functions and their Lie derivatives are termed \emph{exhaustive summary}. By verifying the uniqueness of the exhaustive summary, the structural global identifiability of the model is validated.

The drawback of both these approaches lie on the need to know the number of output derivatives over which the identifiability can be verified. To mitigate some of these issues, iterative approaches have been suggested. The \emph{identifiability tableau} proposed in \cite{balsa-canto_iterative_2010} applied it for Taylor series approach while in \cite{chis_structural_2011}, they were applied for the generating series technique to develop the genSSI MATLAB toolbox.

\subsubsection{Isomorphism based approach}
The isomorphism based approach answer the distinguishability question by analyzing the relationship between state-space realizations. This method exploits the fact that, under certain conditions, indistinguishable state space models have locally isomorphic state spaces. So the identifiability analysis works to show that state isomorphisms must have certain properties within the class of state space systems considered. This helps to parametrize indistinguishable state space models and if the isomorphism can be shown to be identity, then global identifiability is also verified.

One of the earliest works to analyze the identifiability property through the state-space realization theory is in \cite{glover1974parametrizations}. For nonlinear systems, local state isomorphism based works were proposed in \cite{vajda_state_1989} which were followed up by works such as \cite{joly1998some} for uncontrolled systems, \cite{peeters2005identifiability} for homogeneous systems, and \cite{nemcova_structural_2010} for polynomial and rational models with provision for piece-wise continuous inputs. These approaches assume the system model is minimal. However, minimality is not a necessary condition for identifiability. One critique of these approaches is the lack of a systematic method to verify identifiability. While this was mitigated to some extent in the systematic solution proposed in \cite{denis1996identifiability}, this continues to remain an issue to use this method for complex nonlinear models (\cite{chis_structural_2011}). 

\subsubsection{Differential algebraic approach}
The potential of the differential algebraic tools for identifiability problems was discussed in the seminal work in \cite{ljung1994global}. The authors deploy the Ritt's algorithm to find the characteristic set of the polynomial ideal generated by the system model (assuming that it can be written in polynomial form). Given the system model $\Sigma_\theta$ in \eqref{eq_ch_id_gen_nonlinear_model}, the idea is to rewrite the state-space model as a set of polynomials,
\begin{align}
g_i(\frac{d}{dt}, u(t), y(t), \theta) = 0 \label{eq_ch_id_ljung_poly}
\end{align}
along with $\dot{\theta}  = 0$. Here $i = 1, 2,, \cdots$ and $\frac{d}{dt}$ stands for all the higher order derivatives of the inputs and outputs. Then after a careful choice of ranking of the variables, the characteristic set of the ideal generated by the set of polynomials \eqref{eq_ch_id_ljung_poly}  is obtained through Ritt's algorithm. 
The structure of the characteristic set provides an inference of the identifiability characteristic of the original model (local, global, non-identifiable). The authors also show that the identifiability and the estimation of the parameters  are guaranteed when for each parameter $\theta_j$, a linear regression form,
\begin{align}
P_j(\frac{d}{dt}, u(t), y(t)) + \theta_j Q_j(\frac{d}{dt}, u(t), y(t)) = 0
\end{align}
is obtained. The implementation of differential algebra approach has different flavours. One approach, proposed in \cite{saccomani1997global} uses a differential ring that does not consider $\theta$, a strategy framed  in \cite{ollivier1990probleme}. This proves useful for biological systems. Since these models have a large number of parameters, including them in the differential ring incurs significant computational effort. This is elaborated in \cite{audoly2001global} which develops identifiability tools for biological systems.

\subsubsection{Differential geometric approach}
In \cite{tunali1987new}, the authors characterize identifiability as an extended observability problem, where the parameters are added to the state vector and the observability of the new model is evaluated. These results are local in nature, but have an intuitive appeal to it that it lead to the development of the toolbox STRIKE-GOLDD (\cite{villaverde2016structural}).

In \cite{xia2003identifiability}, the implicit function theorem is employed as a means to derive local identifiability results. In particular, local structural identifiability is formulated as algebraic identifiability and illustrated. While the relationship between the various local identifiability characterizations are clearly given, the actual computation steps to validate identifiability is slightly ambiguous and the example provided  seems not well-handled as also noted in \cite{saccomani_effective_2011}.

\subsection{Discrete-time identifiability}
For the discrete-time case, the identifiability results are limited. In \cite{anstett2006chaotic}, the authors formulated the cryptographic key's ability to be cracked as an identifiability problem and then reused the continuous-time results in the discrete-time context. The authors in \cite{anstett_identifiability_2008} developed a discrete-time version of the local state isomorphism theorem and used it to establish identifiability results for discrete-time systems with polynomial nonlinearities.

The authors in \cite{nomm_identifiability_2004} develop discrete-time local identifiability results using the implicit function theorem similar to that of continuous-time in \cite{xia2003identifiability}. These results, however, don't provide any specific systematic procedure, neither do the examples provide any insights into the procedure.

\subsection{Identifiability of LPV models}
For identifiability of LPV models, all the works consider models with static dependence on the scheduling variables. The authors in  \cite{lee1997identifiability} derive some perspectives on those models that could be represented using the linear fractional transformation (LFT) approach. They provide an identifiability characterization of such models using the existence of similarity transformation between two realizations. In \cite{dankers_informative_2011}, the authors deal with the dual problems of identifiability and informativity that concerns the parameter estimation. The models considered are the input-output models with LPV-ARX structure in contrast to the state-space models of interest in this paper. Discrete-time affine LPV model identifiability is discussed in \cite{alkhoury2017identifiability}, where the authors use the realization theory developed in \cite{petreczky2012affine}. The authors provide a systematic procedure that culminates in a rank condition that would verify the presence of an isomorphism between the realizations. The results are necessary and sufficient for local structural identifiability and sufficient for global structural identifiability.

\subsection{Software packages for identifiability evaluation}
While there are several systematic approaches to validate identifiability for different systems, implementation of each of those methods for comparison purposes is difficult. In this respect, the continuous-time results obtained in this paper are compared to that with DAISY (\cite{bellu2007daisy}).

DAISY is a package developed under the REDUCE platform and implements the ideas that originated in \cite{saccomani1997global} and elaborated in \cite{audoly2001global}. The package uses the Ritt's algorithm to eliminate the states of the system and compute the characteristic set associated with the differential ideal generated by the system differential equations. The differential ring used is $\mathbb{R}[x,y,u]$ (instead of the $\mathbb{R}[x,y,u,\theta]$ as in \cite{ljung1994global}) and hence a normalized input-output relation is obtained from the characteristic set. The exhaustive summary \cite{walter1982global} is extracted from the normalized input-output relations by gathering the functions of parameters that appear as coefficients. Further, the authors assign random numerical values to the parameters and subject it to the Buchberger algorithm to compute their Gr\"{o}bner basis (\cite{buchberger2006bruno}). Depending upon the number of solutions it admits, the original system is globally, locally or non-identifiable. 

There are also other packages such as genSSI (\cite{chis_structural_2011}) and STRIKE-GOLDD (\cite{villaverde2016structural}) developed under under the MATLAB computing platform. genSSI evaluates the identifiability of models in the control affine form using the generating series approach from \cite{walter1982global}. And the STRIKE-GOLDD is based on the extended observability approach of \cite{tunali1987new} and computes $n+q$ derivatives of the output and then evaluates the Jacobian of the resulting equation with respect to the extended state vector that includes the unknown parameters.


\subsection{Motivation for the present work}
The above summary of the literature illustrates a wide-range of works that have been carried out for identifiability of nonlinear models. However, the literature is limited when it comes to LPV models. The key points addressed in this paper are:
\begin{itemize}
	\item To develop a procedure to verify identifiability of LPV (and quasi-LPV) models
	\item To explore a unified procedure for a class of both continuous-time and discrete-time LPV models
	\item To utilize and exploit the theoretical and applied results already in the literature
\end{itemize}

With this in mind and the wide-ranging models that can be represented through LPV/quasi-LPV models, the elimination strategy seems appropriate in this context as they can work, to some extent, agnostic to the underlying model. One of the underlying themes in the literature of elimination techniques is to arrive at the exhaustive summary of a model. In \cite{walter1982global}, it is through a generating series, whereas in \cite{audoly2001global}, it is through differential algebra (Ritt's algorithm). In this work, this is achieved using a parity-space based approach.

Before the approach is discussed, the specifications of the problem is discussed in the following section.

\section{Problem formulation} \label{sec_problem_formulation}
In this section, the definitions of identifiability of interest and the assumptions underlying the proposed procedure are given.

\begin{definition}{\textbf{Identifiability}} \label{defn_ch_id_identifiability_general}
A parameter $\theta_i$, $i \in \{1, \cdots, q\}$, is structurally globally (or uniquely) identifiable if for almost any $\theta_i^* \in \Theta$,
\begin{align*}
\Sigma_{\theta_i} = \Sigma_{\theta_i^*} \Rightarrow \theta_i = \theta_i^* .
\end{align*}

A parameter $\theta_i$, $i \in \{1, \cdots, q\}$, is structurally locally identifiable if for almost any $\theta_i^* \in \Theta$, there exists a neighbourhood $\eta(\theta^*)$ such that 
\begin{align*}
\Sigma_{\theta_i} = \Sigma_{\theta_i^*} \Rightarrow \theta_i = \theta_i^* .
\end{align*}

Consequently, if this neighbourhood does not exist, the parameter $\theta_i$ is structurally non-identifiable.
\end{definition}
The above definition formalizes the distinguishability property through structural identifiability. To \emph{a priori} verify identifiability using standard mathematical tools, more tangible definitions are required. In this respect two approaches of interest are discussed below, namely, structural identifiability from \cite{audoly2001global} and algebraic identifiability from \cite{xia2003identifiability}. This characterization requires the following notations and terminologies:
\paragraph{\textbf{Exhaustive Summary} (\cite{walter1982global})} of an experiment is a set of functions,  $\Pi(\theta)$, if it contains only, but all the information about $\theta$ that can be extracted from the knowledge of $u$ and $y$. That is, they embody the parameter dependence of the input-output model completely. These are also referred to as  the observational parameter vector in \cite{jacquez1985numerical}. Some authors use a slightly different terminology, for example, in \cite{audoly2001global}, the authors refer to the set of equations,
	\begin{align}
	{\Pi}({\theta}) = {\Pi}(\tilde{\theta}) \label{eq_ch_id_ex_summ_audoly}
	\end{align}
	as exhaustive summary, where $\tilde{\theta}$ refers to the specific instance of $\theta$ used to verify if $\Pi(\theta)$ admits only one solution, that is, $\tilde{\theta}$. In this work, we use \emph{exhaustive summary} to refer to the generic set of equations denoted by $\Pi(\theta)$ whereas \eqref{eq_ch_id_ex_summ_audoly} would be referred to as  \emph{exhaustive summary evaluation}. 
	
\paragraph{$\Phi(.)$ \textbf{Set of identifiability equations} (\cite{xia2003identifiability})}: These are $q$ equations which are functions of the known and measured variables along with their derivatives, and the unknown parameters. They are of the form:
	\begin{align*}
	\Phi(\theta, y, \dot{y}, \ddot{y}, \cdots, u, \dot{u}, \cdots) = 0
	\end{align*}
	
Given a system model \eqref{eq_ch_id_gen_nonlinear_model}, it is possible to obtain its exhaustive summary through various methods and then validate the number of solutions in $\theta$ admitted by it. This characterization is formalized as follows. Note that the use of $y(\Pi(\theta), t)$ is in reference to the experiment which provides a set of measurement (outputs) that depend on the exhaustive summary.
\begin{definition}[Structural identifiability (\cite{audoly2001global})] \label{defn_ch_id_structural_id}
 A parameter $\theta_i$ is, 
 
  globally (or uniquely) identifiable if and only if, for almost any $\tilde{\theta}$, the following system has only one solution, $\theta_i = \tilde{\theta}_i$, $i \in \{1, \cdots, q\}$:
 \begin{align}
 y(\Pi(\theta), t) = y({\Pi}(\tilde{\theta}),t) \label{eq_ch_id_distinguishability_exhaustive_summary}
 \end{align}
 
  locally (nonuniquely) identifiable (LSI) if and only if, for almost any $\tilde{\theta}$, the system \eqref{eq_ch_id_distinguishability_exhaustive_summary} has (for $\theta_i$) more than one, but a finite number of solutions. \\

  non-identifiable if and only if, for almost any $\tilde{\theta}$, the system \eqref{eq_ch_id_distinguishability_exhaustive_summary} has (for $\theta_i$) infinite number of solutions.  
\end{definition}
The second definition of interest is the algebraic identifiability in \cite{xia2003identifiability}.
\begin{definition}[Algebraic identifiability (AI) (\cite{xia2003identifiability})] \label{defn_ch_id_algebraic_id}
The system model $\Sigma_\theta$ is said to be algebraically identifiable if there exist a $T > 0$, a positive integer $k$ and a meromorphic function
$\Phi : \mathbb{R}^q \times \mathbb{R}^{(k+1)m} \times \mathbb{R}^{(k+1)p} \rightarrow \mathbb{R}^q $ such that \vspace{-5pt} 
\begin{align}
\text{det} \frac{\partial\Phi}{\partial \theta} \neq 0
\end{align}
and 
\begin{align} 
\Phi(\theta, y, \dot{y}, ..., y^{(k)}, u, \dot{u}, ..., u^{(k)}) = 0 \label{eq_ch_id_id_equations}
\end{align}
hold, on $[0,T]$, for all $(\theta,u,... u^{(k)}, y, \dot{y}, ..., y^{(k)})$ where $(\theta, x_0, u)$ belong to an open and dense subset.
\end{definition}

The relationship between algebraic identifiability and the local structural identifiability is clarified in the following proposition. This is done by reiterating the characterization of the two definitions to illustrate the equivalence.
\begin{proposition} \label{prop_ch_id_ai_lsi} For a system model of type $\Sigma_\theta$, the definitions of algebraically identifiable (AI) and locally structurally identifiable (LSI) are equivalent\footnote{Generically, for almost all cases}.
\end{proposition}

To verify this equivalence, the first step would be to consider how one can obtain a set of $q$ equations in the form of $\Phi$. In the procedure described in this paper as well as those in \cite{audoly2001global}, the elimination of states would yield $p$ equations (let's denote them as $\Psi$). If $p<q$, then one could obtain more equations by differentiating $\Psi$ to obtain further equations such that $p=q$. With this set, one can readily verify Proposition 1 by checking that the exhaustive summary satisfies the two definitions (SI and ASI) if only if and the set of identifiability equations also satisfy them. This objective of this proposition would be put to use when we consider equivalent approaches to verify identifiability locally.

\begin{assumption}[Model structure]
The models of interest are those nonlinear models that could be written in the LPV or quasi-LPV form with affine parametrization. That is,
\begin{align}
\dot{x}(t) &= A(\rho(t), \theta) x(t) + B(\rho(t), \theta) u(t) \nonumber \\
\dot{y}(t) &= C(\rho(t), \theta) x(t) + D(\rho(t), \theta) u(t) \label{eq_ch_id_qlpv_model_form}
\end{align}
with $x \in \mathbb{R}^n$, $u \in \mathbb{R}^m$, $y \in \mathbb{R}^p$, $\rho \in \mathbb{R}^{\xi}$ $\theta \in \mathbb{R}^q$  with the appropriate dimensions for the system matrices ($A$, $B$, $C$, $D$) which are of the form, 
\begin{align}
X(\rho(t), \theta) = X_0(\rho(t)) + \sum_{j=1}^{q} \theta_j \bar{X}_j(\rho(t))
\end{align}
\end{assumption}
The scheduling or premise variable, $\rho(t)$, is either composed of external variables with static dependence (in which case the model is LPV) or that of system variables such as an inputs, states and outputs (in which case the model is quasi-LPV). 
\begin{assumption}[Premise variables]
The premise variables of the quasi-LPV model are known or measured.
\end{assumption}

\begin{remark}
The nonlinear models of the form \eqref{eq_ch_id_gen_nonlinear_model} can be rewritten into quasi-LPV forms using several of the existing embedding techniques (see for example, \cite{ohtake2003fuzzy}, \cite{kwiatkowski_automated_2006}, \cite{abbas2014embedding}). The quasi-LPV representation is not unique, and one might obtain models with different types of premise variables. In this work, only those models with known or measured premise variables are considered.
\end{remark}

\begin{assumption}[Characteristics of $f$ and $h$ in \eqref{eq_ch_id_gen_nonlinear_model}]
The state and the output functions, $f$ and $h$ respectively, are assumed to be meromorphic. Further, 
\end{assumption}
\begin{align}
\textrm{rank} \left( \frac{\partial h(x,\theta,u)}{\partial x} \right) = p 
\end{align}
The assumptions on the model functions are a superset to the assumptions given in \cite{audoly2001global} and \cite{xia2003identifiability} to complete the definitions \ref{defn_ch_id_structural_id} and \ref{defn_ch_id_algebraic_id}. In terms of the quasi-LPV model, this condition requires the following
\begin{itemize}
	\item The nonlinearities that appear in the matrices $A(.)$, $B(.)$, $C(.)$, and $D(.)$ are meromorphic.
	\item The rows of the matrix $C(.)x + D(.)u$ are locally independent, that is,
	\begin{align*}
	\text{rank} \left( \frac{\partial}{\partial x} \left[C(.) x + D(.) u\right] \right) = p
	\end{align*}
\end{itemize}

\begin{assumption}[Initial conditions]\label{ass_ch_id_init_conditions}
The state initial conditions are arbitrary.
\end{assumption}

\begin{assumption}[Inputs]
The higher order derivatives of the inputs are defined and are known
\end{assumption}
This assumption is required at least up to the order required for identifiability analysis so that it is possible to formulate $\Phi(.)$ in \eqref{eq_ch_id_id_equations}. 

\begin{remark}[Discrete-time case]
The discussion in this section has focused on continuous-time models, though it holds for the discrete-time case with the exchange of shift in discrete-time for derivatives in continuous-time as commented in \cite{anstett2006chaotic}.
\end{remark}

\section{Parameter identifiability: continuous-time models} \label{sec_id_algo_continuous}
In this section, an overview of the proposed parity-space based identifiability analysis method is given. The method is illustrated with a set of examples and the results obtained are compared with that from DAISY.

\subsection{Steps-by-step description}
The procedure for the identifiability analysis proposed is inspired by the parity-space approach in \cite{chow1984analytical} as a means to eliminate the states of the system. 
The procedure could be summarized as follows:
\paragraph{Step 1: Formulation of algebraic equations} The LPV/quasi-LPV model in \eqref{eq_ch_id_qlpv_model_form} is rewritten for illustration purposes as follows, 
\begin{align}
x^{(1)} &= A^{(0)} x^{(0)} + B^{(0)} u^{(0)} \nonumber \\
y^{(0)} &= C^{(0)} x^{(0)} + D^{(0)} u^{(0)}
\end{align}
where the superscript refers to the order of derivatives, that is,
\begin{align*}
A^{(j)} = \frac{d^{j} \left( A(\rho(t), \theta) \right)}{dt^j} 
\end{align*}
with the dependence on time, premise variables and the parameters are omitted for the sake of simplicity. If the model has to be considered up to $2$nd order derivatives of the output, it is possible to rewrite the above as (with known and measured parts on the left hand side), 
\begin{align}
{\footnotesize{
\begin{bmatrix}
y^{(0)} \\ 
y^{(1)} \\
y^{(2)}  \\
\mathbf{0} \\  
\mathbf{0}\\  
\end{bmatrix}
+
\begin{bmatrix}
- D^{(0)} & \mathbf{0} & \mathbf{0}\\ 
- D^{(1)} & - D^{(0)} & \mathbf{0} \\
- D^{(2)} & - 2 D^{(1)} & - D^{(0)} \\
B^{(0)} & \mathbf{0} & \mathbf{0} \\  
B^{(1)} & \mathbf{0} &  B^{(0)} \\  
\end{bmatrix}
\begin{bmatrix}
u^{(0)} \\
u^{(1)} \\
u^{(2)}
\end{bmatrix}
=
\begin{bmatrix}
C^{(0)} &  \mathbf{0} &  \mathbf{0} \\
C^{(1)} & C^{(0)} &  \mathbf{0} \\
C^{(2)} & 2C^{(1)} & C^{(0)}\\
-A^{(0)} & I_n &  \mathbf{0} \\
-A^{(1)} & -A^{(0)} & I_n \\
\end{bmatrix}
\begin{bmatrix}
x^{(0)} \\ x^{(1)} \\ x^{(2)} 
\end{bmatrix} }} \label{eq_ch_id_algebraic_eqs_form}
\end{align}
More generally for up to an order $w$ of the output derivative,
\begin{align}
\begin{bmatrix}
\mathbb{Y} \\ \mathbf{0}_{w \times n}
\end{bmatrix}
+
\begin{bmatrix}
-\mathbb{D}(\theta) \\ \mathbb{B}(\theta)
\end{bmatrix}
\mathbb{U}
=
\begin{bmatrix}
\mathbb{C}(\theta) \\ \mathbb{A}(\theta)
\end{bmatrix}
\mathbb{X} \label{eq_ch_id_id_form_condensed}
\end{align}
with the left hand side containing known and measured terms. The presence of $(\theta)$ indicates the explicit appearance of the parameter in the matrices. Notice, however that, all the elements except $\mathbb{U}$ are indirectly dependent on the parameter $\theta$. Here,
\begin{align}
\mathbb{Y} &= \begin{bmatrix}
(y^{(0)})^T & (y^{(1)})^T & (y^{(2)})^T & \cdots & (y^{(w)})^T
\end{bmatrix}^T \nonumber \\
\mathbb{U} &= \begin{bmatrix}
(u^{(0)})^T & (u^{(1)})^T & (u^{(2)})^T & \cdots & (u^{(w)})^T
\end{bmatrix}^T \nonumber \\
\mathbb{X} &= \begin{bmatrix}
(x^{(0)})^T & (x^{(1)})^T & (x^{(2)})^T & \cdots & (x^{(w)})^T
\end{bmatrix}^T 
\end{align}
and, 
{\footnotesize{
\begin{align} \label{eq_continuous_matrix_deriv_1}
\mathbb{B}(\theta) &= 
\begin{bmatrix}
B^{(0)} &  \mathbf{0} &  \mathbf{0} &  \mathbf{0} & \cdots &  \mathbf{0}  &\mathbf{0} \\
B^{(1)} & B^{(0)} &  \mathbf{0} &  \mathbf{0} & \cdots &  \mathbf{0}&  \mathbf{0} \\
B^{(2)}  & 2B^{(1)} & B^{(0)}  &  \mathbf{0} & \cdots &  \mathbf{0} &  \mathbf{0}\\
B^{(3)} & 3B^{(2)}  & 3B^{(1)} & B^{(0)}  &  \cdots &  \mathbf{0} &  \mathbf{0} \\
\vdots & \vdots & \vdots & &  \ddots & & \vdots \\
B^{(w-1)} & \binom w2 B^{(w-2)}  & \binom w3  B^{(w-3)} & \binom w4 B^{(w-4)}  &  \cdots &  B^{(0)} &  \mathbf{0}\\
\end{bmatrix} \\
\mathbb{D}(\theta) &= 
\begin{bmatrix}
D^{(0)} &  \mathbf{0} &  \mathbf{0} &  \mathbf{0} & \cdots &  \mathbf{0}  &\mathbf{0} \\
D^{(1)} & D^{(0)} &  \mathbf{0} &  \mathbf{0} & \cdots &  \mathbf{0}&  \mathbf{0} \\
D^{(2)}  & 2D^{(1)} & D^{(0)}  &  \mathbf{0} & \cdots &  \mathbf{0} &  \mathbf{0}\\
D^{(3)} & 3D^{(2)}  & 3D^{(1)} & D^{(0)}  &  \cdots &  \mathbf{0} &  \mathbf{0} \\
\vdots & \vdots & \vdots & &  \ddots & & \vdots \\
D^{(w)} & \binom {w+1}2 D^{(w-1)}  & \binom {w+1}3  D^{(w-2)} & \binom {w+1}4 D^{(w-3)}  &  \cdots & 2 D^{(1)} &  D^{(0)}\\
\end{bmatrix}
\end{align} }}
{\footnotesize{
\begin{align}
\mathbb{C}(\theta) &= 
\begin{bmatrix}
C^{(0)} &  \mathbf{0} &  \mathbf{0} &  \mathbf{0} & \cdots &  \mathbf{0}  &\mathbf{0} \\
C^{(1)} & C^{(0)} &  \mathbf{0} &  \mathbf{0} & \cdots &  \mathbf{0}&  \mathbf{0} \\
C^{(2)}  & 2C^{(1)} & C^{(0)}  &  \mathbf{0} & \cdots &  \mathbf{0} &  \mathbf{0}\\
C^{(3)} & 3C^{(2)}  & 3C^{(1)} & C^{(0)}  &  \cdots &  \mathbf{0} &  \mathbf{0} \\
\vdots & \vdots & \vdots & &  \ddots & & \vdots \\
C^{(w)} & \binom {w+1}2 C^{(w-1)}  & \binom {w+1}3  C^{(w-2)} & \binom {w+1}4 C^{(w-3)}  &  \cdots & 2 C^{(1)} &  C^{(0)}\\
\end{bmatrix} \\
\mathbb{A}(\theta) &= 
\begin{bmatrix}
-A^{(0)} &  I_n &  \mathbf{0} &  \mathbf{0} & \cdots &  \mathbf{0}  &\mathbf{0} \\
-A^{(1)} & -A^{(0)} &  I_n &  \mathbf{0} & \cdots &  \mathbf{0}&  \mathbf{0} \\
-A^{(2)}  & -2A^{(1)} & -A^{(0)}  &  I_n & \cdots &  \mathbf{0} &  \mathbf{0}\\
-A^{(4)} & -3A^{(2)}  & -3A^{(1)} & -A^{(0)}  &  \cdots &  \mathbf{0} &  \mathbf{0} \\
\vdots & \vdots & \vdots & &  \ddots & & \vdots \\
-A^{(w-1)} & -\binom w2 A^{(w-2)}  & -\binom w3 A^{(w-3)} & -\binom w4 A^{(w-4)}  &  \cdots & -A^{(0)} & I_n\\
\end{bmatrix} \label{eq_continuous_matrix_deriv_2}
\end{align} }}
As is apparent, each of these matrices form a Pascal's triangle with the increasing order of derivatives. This is useful when the algorithm is implemented. The representation in \eqref{eq_ch_id_id_form_condensed} is further simplified to indicate the dependence on the unknown parameter as, 
\begin{align}
\mathbb{Y}_0 + \mathbb{G}(\theta) \mathbb{U} = \mathbb{O}(\theta) \mathbb{X} \label{eq_ch_id_id_form_simplified}
\end{align}
Notice that the matrix $\mathbb{O}(\theta)$ has a dimension of $(w p + (w-1)n)\times((w-1)n)$.

\paragraph{Step 2: Computation of null-space}  Once a set of algebraic equations are formulated, the next step is to eliminate the state variables and their derivatives. This is achieved if the left null-space of $\mathbb{O}(\theta)$ is computed, that is, to find a matrix $\Omega(\theta)$, such that,
\begin{align*}
\Omega^T(\theta) \mathbb{O}(\theta) = \mathbf{0}
\end{align*}
For a given $\mathbb{O}(\theta)$, the null-space $\Omega(\theta)$, if it exists, can be computed using symbolic computations such as under the symbolic computation toolbox under MATLAB. The existence of the null-space is directly related to the output and the state matrices which populate $\mathbb{O}(\theta)$.
\paragraph{Step 3: Formulation of the Input-Output-Parameter (I-O-P) equations}
Once the null-space has been obtained, then one can compute from \eqref{eq_ch_id_id_form_simplified},
\begin{align}
\Omega^T(\theta) \left( \mathbb{Y}_0 + \mathbb{G}(\theta) \mathbb{U} \right) = \mathbf{0}
\end{align}
which can alternatively be represented as,
\begin{align}
\Psi(\theta, y , \cdots,  y^{(w)}, u, \cdots, u^{(w)}) = \mathbf{0}
\end{align}
where $\Psi(.)$ is termed as Input-Output-Parameter (I-O-P) equations to signify its dependence on inputs, outputs and parameters, and their derivatives.

\paragraph{Step 4: Identifiability evaluation} Once the I-O-P equations are obtained, the verification of identifiability is achieved through one of the following approaches:
\begin{itemize}
	\item  Following the final step in the DAISY package in \cite{bellu2007daisy}:
	\begin{itemize}
		\item Extract the coefficients of $\Psi(.)$ considering it as polynomials in inputs, outputs, and their derivatives. Those coefficients that depend on the parameters $\theta$ form the exhaustive summary $\Pi(\theta)$
		\item Assign symbolic values for each of the parameters $\{\theta_1, \cdots, \theta_q\}$ and evaluate the exhaustive summary to obtain ${\Pi}(\tilde{\theta})$. For large scale problems, symbolic values can be replaced with numerical values. 
		\item Apply Buchberger algorithm on ${\Pi}(\tilde{\theta})$ to obtain all the solutions. Depending upon the number of solutions ${\Pi}(\tilde{\theta})$ admits, identifiability can be evaluated using Definition \ref{defn_ch_id_structural_id}. In case of numerical approach, the last two steps are repeated several times. Since the results are \emph{generic}, that is, valid for almost all numerical values except for a set of measure zero. This repetition would help to avoid reaching conclusions based on possibly choosing a numerical value of this set of measure zero.
	\end{itemize}
	\item A local identifiability verification using Jacobians (inspired by Proposition 1)
	\begin{itemize}
		\item Obtain the set of $q$ Identifiability equations $\Phi(.)$. Note that the number of I-O-P equations is equal to the number of outputs $p$ and so if $p<q$, one has to differentiate the I-O-P equations and set $\dot{\theta}=0$ to obtain $q$ equations (i.e., $\Phi$).
		\item  	Compute the Jacobian of $\Psi(.)$ with respect to the parameters $\theta$, that is,
	\begin{align*}
	\text{rank} \left( \frac{\partial \Psi}{\partial \theta} \right) = q
	\end{align*}
		\item If the rank is $q$, then local identifiability is verified.
	\end{itemize}
\end{itemize}

\subsection{Illustrative examples}
In this section, several examples are given to show the steps involved in the parity-space approach. The results obtained from the parity-space approach are validated by comparing it with that of DAISY, STRIKE-GOLDD and genSSI software packages. Further, the intermediate step of exhaustive summary is compared with that obtained using DAISY.
\begin{example} This example is used to show all the steps of the proposed approach. Further, it also concerns a model which is not identifiable. Consider the second order nonlinear model,
\begin{align*}
\dot{x}_1 &=\theta_1 x_1 + \theta_2 x_2 u, \\
\dot{x}_2 &= \theta_3 x_1 - x_2, \\
y &= x_1 u,
\end{align*}
\end{example}
A quasi-LPV equivalent form with $x = [x_1\ \ x_2]^T$ and $\rho = u$ is,
\begin{align*}
\dot{x} &= 
\begin{bmatrix} \theta_1 & \theta_2 u \\ \theta_3 & -1\end{bmatrix}
x \qquad \text{and} \qquad
y = \begin{bmatrix}
u & 0
\end{bmatrix}
x
\end{align*}
Based on the specifications of the model structure required, the genSSI software cannot handle this example. Using the parity-space approach with output up to $\ddot{y}$, we obtain the following representation to that in \eqref{eq_ch_id_id_form_simplified},
\begin{align*}
\begin{bmatrix}
 y \\ \dot{y} \\  \ddot{y} \\   0 \\   0 \\   0 \\   0
\end{bmatrix}
=
\begin{bmatrix}
u & 0 & 0 & 0 & 0 & 0 \\
\dot{u} & 0 & u & 0 & 0 & 0 \\
\ddot{u} & 0 & 2 \dot{u} & 0 & u & 0 \\
-\theta_1 & -\theta_2 {u} & 1 & 0 & 0 & 0 \\
-\theta_2 & -\theta_3 & 0 & 1 & 0 & 0 \\
0 & -\theta_2 \dot{u} & -\theta_1 & -\theta_2 u & 1 & 0 \\
0 & 0 & -\theta_2 & \theta_3 & 0 & 1
\end{bmatrix}
\begin{bmatrix}
x_1 \\ x_2 \\ \dot{x}_1 \\ \dot{x}_2 \\ \ddot{x}_1 \\ \ddot{x}_2
\end{bmatrix}
\end{align*}
The left null-space of the matrix $\mathbb{O}(\theta)$ is given by,
\begin{align*}
\Omega^T(\theta) = 
\begin{bmatrix}
(\theta_2^2 u^3 + u \ddot{u}  - 3 \dot{u}^2 + \theta_1  \theta_3 u^2  - 2 \theta_1 u \dot{u} + \theta_3 u \dot{u})/u^3 \\
(3\dot{u} + \theta_1 u - \theta_3 u)/u^2 \\
-1/u \\
-(\dot{u}-\theta_3 u )/u\\
\theta_2 u \\
1 \\ 
0
\end{bmatrix}^T
\end{align*}
which leads to the I-O-P equation as follows,
\begin{align*}
\Psi (.)= \theta_1 u^2 y - 3 \dot{u}^2 y - u^2 \ddot{y} - u^2 \dot{y} + \theta_1 u^2 \dot{y} + u \dot{u} y + 3 u \dot{u} \dot{y} + u \ddot{u} y - 2 \theta_1 u \dot{u} y + \theta_2 \theta_3 u^3 y
\end{align*}
And the exhaustive summary obtained by extracting the coefficients considering $\Psi(.)$ as a polynomial in inputs, outputs and their derivatives. Considering only those coefficients that depend on various inputs, outputs, their derivatives and their product combinations, the following is obtained,
\begin{align}
\Pi(\theta) = \{1 - 2 \theta_1, \   \theta_1 - 1,\   \theta_1,\   \theta_2 \theta_3\} \label{eq_ch_id_eg1_ex_summary}
\end{align}
To verify the number of solutions admitted by this exhaustive summary, a Gr\"{o}bner basis analysis is performed. One strategy is to assign symbolic values for for each of the parameter ($\tilde{\theta}_1 = a$, $\tilde{\theta}_2 = b$, $\tilde{\theta}_3 = c$) and evaluate the exhaustive summary to obtain the specific exhaustive summary,
\begin{align*}
\{2\theta_1-2a,\  \theta_1-a,\   \theta_1-a,\   \theta_2 \theta_3-bc\}
\end{align*}
For this simple example, it is straightforward to see that only $\theta_1$ is identifiable as it admits a unique solution and the other two parameters $\theta_2$ and $\theta_3$ can have several solutions. Hence the model is not identifiable. To formally verify this, these polynomial equations were given as input to the Buchberger algorithm implemented in MuPAD CAS under MATLAB. The Gr\"{o}bner basis for this set is given by,
\begin{align*}
\{\theta_1-a,\ \theta_2 \theta_3-bc\}
\end{align*}
which, if the equations admit a unique solution should have returned $\theta_i = \tilde{\theta}_i$ for $i = 1,2,3$. However, it is not the case here and so the model is not identifiable (or only $\theta_1$ is identifiable).

\paragraph{Comparison with DAISY} The normalized input-output equation obtained through the DAISY package is given by,
\begin{align*}
\Gamma = \ddot{u} u^4 y - 3 \dot{u}^2 u^3 y + 3 \dot{u} \dot{y} u^4 + \dot{u} u^4 y (-2 \theta_1+ 1) - \ddot{y} u^5 + \dot{y} u^5 (\theta_1-1) + u^6 y \theta_2 \theta_3 + u^5 y \theta_1
\end{align*}
which has the same set of exhaustive summary as given in \eqref{eq_ch_id_eg1_ex_summary}. The DAISY package results also verify those inferred above.
\begin{example} In this example, the case of local identifiability is illustrated.
\begin{align*}
\dot{x} &= 
\begin{bmatrix}
\theta_1 & \theta_2 u \\
\theta_2 & -\theta_3
\end{bmatrix} x \\
y & = \begin{bmatrix}
u & 0
\end{bmatrix} x
\end{align*}
\end{example}
Using parity-space approach, the I-O-P equation obtained  is
\begin{align*}
\Psi(.) = \theta_1 u^2 \dot{y} - u^2 \ddot{y} -3 \dot{u}^2 y - \theta_3 u^2 \dot{y} + \theta_2^2 u^3 y + 3 u \dot{u} \dot{y} + u \ddot{u} y - 2 \theta_1 u \dot{u} y + \theta_3 u \dot{u} y + \theta_1 \theta_3 u^2 y
\end{align*}
which has the exhaustive summary of,
\begin{align*}
\Pi(\theta) = \{\theta_3-2 \theta_1, \ \theta_1-\theta_3, \ \theta_1 \theta_3, \ \theta_2^2\}
\end{align*}
The Gr\"{o}bner for this summary with the symbolic assignment of $\tilde{\theta}_1 =a$, $\tilde{\theta}_2 =b$ and $\tilde{\theta}_3=c$ was obtained as,
\begin{align*}
\{\theta_1-a, \ \theta_3-c, \ \theta_2^2-b^2\}
\end{align*}
which indicates that while $\theta_1$ and $\theta_3$ are identifiable, $\theta_2$ is only locally identifiable. The results and the exhaustive summary compares with that obtained from DAISY.

\begin{example}[Air Handling Unit] Consider a simple model of a heat exchanger in \cite{srinivasarengan2016takagi}. The model has been simplified by considering that the inlet air temperature and water temperature are known and constant. A quasi-LPV representation of that model is given by,
\begin{align*}
\begin{bmatrix}
\dot{x}_1  \\ \dot{x}_2
\end{bmatrix}
&=
\begin{bmatrix}
 -\theta_1 u_1 -\theta_2  & \theta_2 \\
 \theta_4 & -\theta_3 u_2 - \theta_4
\end{bmatrix}
\begin{bmatrix}
x_1 \\ x_2
\end{bmatrix}
+
\begin{bmatrix}
\theta_1 & 0 \\ 0 & 5 \theta_3
\end{bmatrix}
\begin{bmatrix}
u_1 \\ u_2
\end{bmatrix} \\
y &= \begin{bmatrix}
1 & 0
\end{bmatrix}
\begin{bmatrix}
x_1 \\ x_2
\end{bmatrix}
\end{align*}
\end{example}
The exhaustive summary obtained from the parity space approach for this model is
\begin{align*}
\Pi(\theta) &= \{3-\theta_4-\theta_2,\  \theta_2 + \theta_4 - 2,\  \theta_1,\  -\theta_1,\  5\theta_2\theta_3,\  -\theta_3,\ \\
 & \qquad \qquad \qquad \qquad   \theta_3 - \theta_2\theta_3,\  \theta_1\theta_4,-\theta_1,\  2\theta_1 - \theta_1\theta_4,\   \theta_1\theta_3,\  -\theta_1\theta_3 \}
\end{align*}
For large scale models, performing Gr\"{o}bner analysis with symbolic values for the parameters could become intractable. In such cases, especially when practical applications are involved (where the range over which the parameters can take value is predictable), one can reliably use numerical values. By choosing arbitrary numerical values, $\tilde{\theta}_1 = 1$, $\tilde{\theta}_2 = 2$, $\tilde{\theta}_3 = 3$, $\tilde{\theta}_4 = 5$, the specific instance of the exhaustive summary was obtained and the Gr\"{o}bner basis obtained is
\begin{align*}
\{ \theta_2-2, \ \theta_4-5, \ \theta_3-3, \ \theta_1-1 \}
\end{align*}
indicating that the model is globally identifiable. And these results verify with those obtained by DAISY both for the exhaustive summary and the eventual identifiability interpretation. Because a numerical value was used, the results are local in nature. Further, as suggested in \cite{bellu2007daisy}, to obtain confidence on the obtained results, this analysis should be repeated for several sets of arbitrary numerical values.

\section{Parameter identifiability for discrete-time models} \label{sec_id_algo_discrete}
In practical scenarios, parameter estimation involves discrete-time models. Hence it is vital to consider the identifiability of system models in discrete-time. In this section, a brief outline to extend the parity-space method to discrete-time quasi-LPV models is given. It is to be noted that the effect of discretization on the identifiability is not treated here. Consider a discrete-time quasi-LPV model of the form,
\begin{align}
x_{k+1} &= A(\rho_k, \theta) x_k + B(\rho_k, \theta)\nonumber \\
y_{k} &= C(\rho_k, \theta) u_k + D(\rho_k, \theta) \label{eq_ch_id_discrete-time_quasi_LPV}
\end{align}
For this type of model, the procedure for parameter identifiability follows a similar trajectory. The key difference is in the first step where the set of algebraic equations is obtained in a different way. For the sake of simplicity, in the following, $A_k$ would be used in place of $A(\rho_k,\theta)$ and similarly for other matrices. \\

For the discrete-time case, the algebraic equations take a far simpler structure compared to that in the continuous-time case. The continuous-time algebraic equations in \eqref{eq_ch_id_id_form_condensed} is rewritten for the discrete-time case as:
\begin{align} 
\begin{bmatrix}
\mathbb{Y}_k \\ \mathbf{0}_{w \times n}
\end{bmatrix}
+
\begin{bmatrix}
-\mathbb{D}_k \\ \mathbb{B}_k
\end{bmatrix}
\mathbb{U}_k
=
\begin{bmatrix}
\mathbb{C}_k \\ \mathbb{A}_k
\end{bmatrix}
\mathbb{X}_k \label{eq_ch_id_id_form_condensed_discrete-time}
\end{align}
where,
\begin{align}
\mathbb{Y}_k &= \begin{bmatrix} y_{k}^T & y_{k+1}^T & \cdots & y_{k+w}^T \end{bmatrix}^T \nonumber \\
\mathbb{U}_k &= \begin{bmatrix} u_{k}^T & u_{k+1}^T & \cdots & u_{k+w}^T \end{bmatrix}^T \nonumber \\
\mathbb{X}_k &= \begin{bmatrix} x_{k}^T & x_{k+1}^T & \cdots & x_{k+w}^T \end{bmatrix}^T
\end{align}
and

\begin{align} \label{eq_discrete_matrix_deriv_1}
\mathbb{B}_k &= 
\begin{bmatrix}
B_{k} &  \mathbf{0}  &  \mathbf{0} & \cdots &  \mathbf{0}  &\mathbf{0} \\
\mathbf{0} & B_{k+1}  &  \mathbf{0} & \cdots &  \mathbf{0}&  \mathbf{0} \\
\mathbf{0} &   \mathbf{0}  & B_{k+2} &   \cdots &  \mathbf{0} &  \mathbf{0}\\
\vdots & \vdots  & &  \ddots & & \vdots \\
\mathbf{0} & \mathbf{0}  & \mathbf{0}   &  \cdots &  B_{k+w-1} &  \mathbf{0}\\
\end{bmatrix} 
\end{align}
\begin{align}
\mathbb{D}_k & = 
\begin{bmatrix}
D_{k} &  \mathbf{0}  &  \mathbf{0} & \cdots &  \mathbf{0}  &\mathbf{0} \\
\mathbf{0} & D_{k+1}  &  \mathbf{0} & \cdots &  \mathbf{0}&  \mathbf{0} \\
\mathbf{0} &   \mathbf{0}  & D_{k+2} &  \cdots &  \mathbf{0} &  \mathbf{0}\\
\vdots & \vdots  & &  \ddots & & \vdots \\
\mathbf{0} & \mathbf{0}  & \mathbf{0}   &  \cdots  &  \mathbf{0} &  D_{k+w}\\
\end{bmatrix} 
\end{align}
\begin{align}
\mathbb{A}_k &= 
\begin{bmatrix}
-A_{k} &  I_n  &  \mathbf{0} &  & \mathbf{0}\cdots &  \mathbf{0}  &\mathbf{0} \\
\mathbf{0} & -A_{k+1}  &  I_n & & \mathbf{0}\cdots &  \mathbf{0}&  \mathbf{0} \\
\mathbf{0} &   \mathbf{0}  & -A_{k+2} & I_n   & \cdots &  \mathbf{0} &  \mathbf{0}\\
\vdots & \vdots  & \vdots & &  \ddots & & \vdots \\
\mathbf{0} & \mathbf{0}  & \mathbf{0}  & \mathbf{0}  &  \cdots &  A_{k+w-1} &  I_n\\
\end{bmatrix}
\end{align}
\begin{align} 
\mathbb{C}_k &= 
\begin{bmatrix}
C_{k} &  \mathbf{0}  &  \mathbf{0} & \cdots &  \mathbf{0}  &\mathbf{0} \\
\mathbf{0} & C_{k+1}  &  \mathbf{0} & \cdots &  \mathbf{0}&  \mathbf{0} \\
\mathbf{0} &   \mathbf{0}  & C_{k+2}   & \cdots &  \mathbf{0} &  \mathbf{0}\\
\vdots & \vdots  & &  \ddots & & \vdots \\
\mathbf{0} & \mathbf{0}  & \mathbf{0}   &  \cdots  &  \mathbf{0} &  C_{k+w}\\
\end{bmatrix}  \label{eq_discrete_matrix_deriv_2}
\end{align}

The other three steps in this case follow that of the continuous time approach with the derivatives replaced with time shifts (Steps 2 to 4 in Sec.\ref{sec_id_algo_continuous}).

\subsection*{Illustrative examples}
All the software packages in the literature are available only for continuous-time models. Hence, comparison of results with any existing packages is not feasible. Hence, the examples are picked from existing literature on discrete-time identifiability and the results compared with those obtained from methods proposed in this work.
\begin{example}This is an example of the Henon map adapted from \cite{anstett2006chaotic} 
\begin{align*}
x_{1,k+1} &= \theta_1 x_{1,k}^2 + \theta_2 x_{2,k} + u_k \\
x_{2,k+1} &= \theta_3 x_{1,k} + \theta_4 u_k \\
y_k &= x_{1,k} 
\end{align*}
\end{example}
The I-O-P equation obtained is:
\begin{align*}
\Psi(.) = -\theta_2 \theta_3 y_k^2 - u_{k+1}  + y_{k+2} -& u_{k}(\theta_2 \theta_4 + \theta_1 y_{k+1})  + \theta_1 y_{k+1} (u_{k} - y_{k+1})
\end{align*}
The exhaustive summary obtained using the parity-space approach is:
\begin{align*}
\Pi(\theta) = \{\theta_1, \ \theta_2 \theta_4, \ \theta_2 \theta_3\}
\end{align*}
The identifiability results verify with that from \cite{anstett2006chaotic} that only the parameter $\theta_1$ is identifiable.
\begin{example}This is also an example from \cite{anstett2006chaotic} of Burgers map
\begin{align*}
x_{1,k+1} &= (1+\theta_1) x_{1,k} + x_{1,k} x_{2,k} + u_k \\
x_{2,k+1} &= (1-\theta_2) x_{2,k} - x_{1,k}^2 u_k  \\
y_k &= x_{1,k}
\end{align*}
\end{example}
The I-O-P equation obtained for this case is:
\begin{align*}
\Psi(.) = y_{k+2} &- u_{k+1} + (y_{k+1} - \theta_2 y_{k+1})(\theta_1 + 1) - u_{k}(1- y_{k+1} y_k^2 + \theta_1) \\ +& \frac{(u_k - y_{k+1}) (y_{k} + y_{k+1} + \theta_1 y_{k} - \theta_2 y_{k+1})}{y_k}
\end{align*}
The exhaustive summary obtained for this example is:
\begin{align*}
\Pi(\theta) = \{\theta_1, \theta_2, \theta_1-\theta_2-\theta_1\theta_2\}
\end{align*}
It is easy to see that the model is identifiable and agrees with the results in \cite{anstett2006chaotic}.

\section{Towards a systematic formulation} \label{sec_ch_id_algo}
In this section, we discuss how to realize a systematic implementation of the proposed algorithm. This includes some algorithmic steps for sample scenarios. The implementation and realisation as a toolbox is for a future work. These formulations are for both continuous-time and discrete-time models with appropriate modifications, though the discussion focuses on continuous-time models.

\subsection{The choice on the number of derivatives} The discussion in the preceding sections did not explicitly consider $w$,  the number of derivatives (or shifts in discrete-time) for which the null-space $\Omega^T(\theta)$ exists and  hence the I-O-P equations and the exhaustive summary that follow. This corresponds to the observability index of the system model. A detailed discussion on observability index of a nonlinear system could be referred to in \cite{nijmeijer1990nonlinear} though a brief idea is given below. Consider a SISO system of the form \eqref{eq_ch_id_gen_nonlinear_model} and the observability index of this model is defined as $w>0$, if in the neighbourhood of $x_0$,
\begin{align*}
\text{rank} \left[L_f^{w-1} h, L_f^{w-2} h, \cdots, L_f^{0} h \right] = w \\
\text{and} \qquad \qquad\text{rank} \left[L_f^{w} h, L_f^{w-1} h, \cdots, L_f^{0} h \right] = w
\end{align*}
where, $L_f h$ corresponds to the Lie derivative of $h$ over $f$, that is,
\begin{align*}
L_f h \triangleq \frac{\partial h(x,u,\theta)}{\partial x} f(x,u,\theta)
\end{align*}
and $L_f^i h$ refers to $i$th successive application of the Lie derivative. Essentially it means that, locally, the dimension of the space spanned by the model does not grow after $(w-1)$ derivatives. For MIMO systems, the observability index is defined for each output. A discrete-time version of this is briefly discussed in Chapter 5 of \cite{anstett2006systemes}.

This means that for a SISO system, using $w$ derivatives would guarantee that $\Omega^T(\theta)$ exists and hence would provide the I-O-P equation corresponding to the output. The next question is to know whether the observability index has been connected to the system dimensions theoretically. That is, is it possible to obtain the index without verifying the rank condition. For linear systems, this is less than or equal tos the number of states (easily to verified using Cayley-Hamilton theorem). For nonlinear systems, for the following two classes:
\begin{itemize}
	\item models of the form \eqref{eq_ch_id_gen_nonlinear_model} where the functions are rational 
	\item models in a control affine form \eqref{eq_ch_id_control_affine} with analytical functions
\end{itemize}
it has been shown that, locally, the observability index has an upper bound, equal to the number of states in the system $n$. That is $n$ derivatives of outputs are sufficient to guarantee that the null-space $\Omega^T(\theta)$ exists. For a more detailed discussion, see \cite{anguelova2007observability}. This is an upper bound because: one, the model may not be minimal and has unobservable spaces and hence the observability index is less than $n$. Second, for MIMO systems, each output would have different observability indices and hence the total number of derivatives required to span the entire observability space can be less than $n$.

Consequently, for single output systems, $n$ derivatives of output would guarantee that the null-space $\Omega^T(\theta)$ exists. Hence $w=n$ for SISO systems. For MIMO systems, this is further complicated. Each output's observability index has an upper bound of $n$, but is more likely to be lower than $n$. A systematic approach to handle this scenario is discussed later in this section.


\subsection{Algorithm for parameter identifiability}
The algorithm that was used for the analysis of identifiability for the illustrative examples discussed in the previous section is summarized in Algo.~\ref{alg_ch_id_identifiability_parity_space}. The implementation was done on the MATLAB computing environment with the use of symbolic computation toolbox and MuPAD computer algebra systems (CAS). Once the set of I-O-P equations is obtained and the exhaustive summary extracted, the Gr\"{o}bner basis evaluation is performed through the MuPAD CAS scripts. Hence, at this moment, there are components of the algorithms that require manual intervention. The first step in the algorithm chooses the upper bound on the number of derivatives to be $n$. While this is applicable for the models chosen for illustrative example, this is not generic other than for the models specified before. Further, for MIMO systems, since the observability index depend on individual outputs, a step by step analysis starting from $0$ derivatives is considered. Further optimization is envisaged in this respect.

\begin{algorithm}
\begin{algorithmic}[1]
\STATE Choose the upper bound on the number of derivatives $w=n$, the matrices $\mathbb{A}, \mathbb{B}, \mathbb{C}, \mathbb{D}$.
\STATE Evaluate the matrices (\eqref{eq_continuous_matrix_deriv_1}-\eqref{eq_continuous_matrix_deriv_2}) and their higher order (element-wise) derivatives 
\FOR{w = 0 to n}
	\STATE Formulate $\mathbb{Y}_0 + \mathbb{G}(\theta) \mathbb{U} = \mathbb{O}(\theta) \mathbb{X}$ as in \eqref{eq_ch_id_id_form_simplified}
	\STATE Compute  $\Omega^T(\theta)$, the left null-space of $\mathbb{O}(\theta)$ using symbolic computation
	\STATE Obtain the I-O-P equations $\Psi(.)$ and extract the coefficients to obtain the exhaustive summary $\Pi(\theta)$.
	\FOR{j = 1 to Number of Iterations}
		\STATE Choose random values for the parameters $\theta_1, \cdots, \theta_q$
		\STATE Evaluate Gr\"{o}bner basis and verify the number of solutions admitted by the exhaustive summary
	\ENDFOR
	\IF{Global or Local identifiability satisfied}
		\STATE END
	\ENDIF
\ENDFOR
\end{algorithmic}
\caption{An algorithm for parameter identifiability}
\label{alg_ch_id_identifiability_parity_space}
\end{algorithm}

\paragraph{Analyzing outputs independently}
One of the assumptions that is part of the problem specifications (and adopted from \cite{xia2003identifiability}) is, 
\begin{align*}
\text{rank} \left( \frac{\partial h(x,u,\theta)}{\partial x} \right) = p
\end{align*}%
That is, the outputs are at least locally independent. This provides an opening to develop local structural identifiability analysis methods that can provide the following advantages:
\begin{itemize}
	\item Obtain the local identifiability results through Jacobian analysis instead of the Buchberger algorithm to obtain the Gr\"{o}bner basis.
	\item As suggested in \cite{bellu2007daisy}, there are $p$ normalized input-output equations. By considering one input at a time, the exit criterion for the algorithm could be set as one I-O-P equation per output by considering the system with one output at a time.
\end{itemize}
Note that, this approach makes sense only if $p \geq q$, as discussed in the remark following Proposition 1. If not, one has to differentiate the I-O-P equations to obtain an appropriate set of Identifiability equations $\Phi(.)$, which involves combinatorial possibilities and is beyond the scope of the current attempt. This limited case is realized as an algorithm in Algo.~\ref{alg_ch_id_local_identifiability_parity_space}.
\begin{algorithm}[ht]
\begin{algorithmic}[1]
\STATE Choose the maximum value for observability index for each output ($w_1, \cdots, w_p$) is $n$.
\STATE Evaluate the matrices (\eqref{eq_discrete_matrix_deriv_1}-\eqref{eq_discrete_matrix_deriv_2}) and their higher order derivatives (element-wise)
\FOR{$i = 1\ \text{to}\ p$ (for each output)}
	\FOR{$w = 0\ \text{to}\ n$}
		\STATE Formulate $\mathbb{Y}_0 + \mathbb{G}(\theta) \mathbb{U} = \mathbb{O}(\theta) \mathbb{X}$ as in \eqref{eq_ch_id_id_form_simplified}
		\STATE Compute $\Omega^T(\theta)$, the left null-space of $\mathbb{O}(\theta)$ using symbolic computation
		\STATE Obtain the I-O-P equation $\psi_i(.)$
		\IF{ one I-O-P equation is obtained}
			\STATE End of search for output $i$
		\ENDIF
	\ENDFOR
	\STATE Add the I-O-P for output to the overall I-O-P, $\Psi(.) = \{\Psi(.), \psi_i(.)\}$
\ENDFOR
\STATE Evaluate $\frac{\partial \Psi}{\partial \theta}$ and compute the rank
\IF{rank $\frac{\partial \Psi}{\partial \theta}$ = q}
	\STATE Model is Locally structurally identifiable
\ELSE
	\STATE Model is not Identifiable
\ENDIF
\end{algorithmic}
\caption{An algorithm for local structural identifiability}
\label{alg_ch_id_local_identifiability_parity_space}
\end{algorithm}

\section{Concluding remarks and perspectives} \label{sec_conclusions}
In this paper, we proposed a procedure for verifying identifiability of LPV and quasi-LPV models in continuous-time and discrete-time models. The procedure exploits the parity-space approach to eliminate the states and uses the residual set of input-output-parameter equations to verify the identifiability of the model. Through several examples the procedures were illustrated and compared with the results obtained from several existing literature. With these preliminary understanding, the algorithm looks to be a useful candidate in the domain LPV model analysis. Given the nature of this paper, there are several avenues for improvements to realize a robust identifiability procedure.

\paragraph{Extending results to newer class of system models} It was noted that the parity-space approach would work well for polynomial parametrization as well. This should be formally extended. While the initial conditions were considered arbitrary in this paper, there are cases where such assumption can be detrimental (See Example 4 in \cite{villaverde2017structural} and Example 2 in \cite{denis-vidal_results_1999}). These cases need to be carefully handled in the implementation to cover a larger spectrum of models and initial conditions. Further, the case of known initial conditions and partially known initial conditions shall also be handled in this extension. The results were restricted to measured or known premise variables. This could be extended to LPV models which has unmeasured or estimated premise variables. 

\paragraph{Systematic implementation} A systematic implementation of the procedure is an immediate future work. The realization could be in the form of a toolbox in MATLAB similar to those such as \cite{chis_structural_2011} or \cite{villaverde2016structural}. That is, a detailed strategy for the model input, options for evaluating local/global identifiability results, and the final display of the results and other relevant information. 

\paragraph{A numerical approach} The STRIKE-GOLDD toolbox \cite{villaverde2016structural} offers to evaluate identifiability either numerically or symbolically. In the numerical approach, a random set of initial conditions are chosen for the states and random values are associated for inputs and their derivatives (these random values are chosen as prime numbers to avoid undesirable cancellations). This significantly reduces the computational effort required to compute the Jacobian. It is to be noted that numerical approach here is not completely numerical. It still requires computing symbolically the Lie derivatives to set up the Observability-Identifiability matrix. A similar inclusion of numerical methods can reduce the symbolic computations. In the Algo.~\ref{alg_ch_id_local_identifiability_parity_space}, this would change the initial assignments step and the Jacobian computation step. The null-space computation steps could also be replaced by exploiting the works in polynomial null-space computation (such as \cite{anaya2009improved} and \cite{khare2010algorithm}).

\paragraph{Computational complexity and efficiency}  The computational complexity analysis of the proposed algorithm would be a topic of future interest once the implementation is realized completely in a computational environment like MATLAB. This would include analysing the relative efficiencies of deploying Buchberger algorithm for Gr\"{o}bner basis computation versus the Jacobian evaluation for local structural identifiability.

Another related interest is the computational comparison with other methods. The DAISY software package envisioned to reduce the computational complexity of the approach in \cite{ljung1994global} by a choice of differential ring that does not contain $\theta$. This works well for biological systems with a number of parameters and a small number of states. However, in engineering systems one often encounters a model with large number of states and a relatively few parameters. An application such as DAISY suffers from the same type of computational overhead as \cite{ljung1994global} had for biological systems. It is to be analyzed whether the parity-space approach can bring in any specific advantages. Similarly to analyze the same type of models for same characteristics, the parity-space approach should also be compared with genSSI and STRIKE-GOLDD toolboxes. Further, as in the recent work by \cite{joubert2020structural}, a numerical sensitivity analysis of the work would be imperative to better understand the numerical properties of the algorithm.

\section*{Acknowledgements}
The authors would like to thank Florian Anstett-Collin for her inputs and discussions.

\bibliography{./../../researchRef}
\bibliographystyle{elsarticle-harv}

\addtolength{\textheight}{-12cm}   



%
%

\end{document}